\documentclass[12pt.]{amsart}
\usepackage{amssymb}
\usepackage{amsthm}
\usepackage{mathtools}
\usepackage{lmodern}
\usepackage[T1]{fontenc}
\usepackage[headings]{fullpage}
\usepackage{paralist}
\usepackage[usenames]{color}
\usepackage{xspace}
\usepackage[normalem]{ulem}
\usepackage[shortalphabetic,abbrev,lite]{amsrefs}
\usepackage{enumerate}

\newtheorem{Theorem}{Theorem}

\newtheorem{Proposition}[Theorem]{Proposition}
\newtheorem{Remark}[Theorem]{Remark}

\newtheorem{Definition}[Theorem]{Definition}

\newtheorem*{Proposition*}{Proposition}

\begin{document}

\title{Betti numbers of piecewiselex ideals}
\author{Christina Jamroz and Gabriel Sosa}

\begin{abstract}
We extend a result of Caviglia and Sbarra to a polynomial ring with base field of any characteristic.  Given a homogeneous ideal containing both a piecewise lex ideal and an ideal generated by powers of the variables, we find a lex ideal with the following property: the ideal in the polynomial ring generated by the piecewise lex ideal, the ideal of powers, and the lex ideal has the same Hilbert function and Betti numbers at least as large as those of the original ideal.
\end{abstract}

\maketitle

\section{Introduction}

Hilbert functions and graded Betti numbers are widely studied invariants in commutative algebra.  In particular, the problem of transforming an ideal into another that has the same Hilbert function and graded Betti numbers greater than or equal to those of the original ideal is one of interest to many researchers.  One of the earliest results in this direction was Macaulay's Theorem \cite{Ma}, which states that if $A=K[x_1,...,x_n]$ is a polynomial ring over a field $K$, then there exists a lex segment ideal realizing the Hilbert function of any homogeneous ideal of $A$.  Later, Bigatti \cite{B}, Hulett \cite{H}, and Pardue \cite{P} proved that lex segments ideals attain the highest Betti numbers among all ideals having the same Hilbert function.

There are two main conjectures in this area of research.  The first is a conjecture of Eisenbud, Green, and Harris \cite{EGH1}, \cite{EGH2}, which asserts that for a homogeneous ideal $I$ containing a homogeneous regular sequence $(f_1,...,f_r)$ with degrees $e_1\leq e_2\leq \cdots \leq e_r$, there exists a lex-plus-powers ideal $L+P$ which has the same Hilbert function as $I$, where $P=(x_1^{e_1},...,x_r^{e_r})$.  The second is Evans' Lex-Plus-Powers Conjecture \cite{FR}, which proposes that in this situation, the graded Betti numbers are such that $b_{ij}(L+P)\geq b_{ij}(I)$ for all $i,j$.  Recently, many people have proven a series of results related to these conjectures, for example \cite{A}, \cite{CK}, \cite{CM}, \cite{C}, \cite{F}, \cite{MM}, \cite{MP}, and \cite{R}.  A strong result has been shown by Mermin and Murai in \cite[Theorem 8.1]{MM}.  They prove the Lex-Plus-Powers Conjecture holds when $(f_1,...,f_r)$ is a regular sequence of monomials.  Notice that, under this assumption, the Eisenbud-Green-Harris Conjecture easily follows from Clements and Lindstr{\"o}m's Theorem \cite{CL}.

A generalization of the Mermin-Murai result was shown by Caviglia and Sbarra \cite{CS}.  In their article, the authors study homogeneous ideals $I$ containing $P+\widetilde{L}$, where $\widetilde{L}$ is a piecewise lex ideal, that is, an ideal which is the sum of extensions to $A$ of lex segment ideals $L_i\subset K[X_1,...,X_i]$.  The quotient rings $A/(P+\widetilde{L})$ are known as Shakin rings.  Their result states that there is a lex ideal $L$ such that $P+\widetilde{L}+L$ and $I$ have the same Hilbert function, and the graded Betti numbers do not decrease when we replace $I$ by the ideal $P+\widetilde{L}+L$.  Unfortunately, the upper bound for the graded Betti numbers was only shown when char($K$)=0.

The main theorem of this paper removes the assumption on the characteristic of the field $K$ in the above result.  In Section 2, we describe the operations done in \cite{MM} to replace the ideal $I$ by a strongly-stable-plus-P ideal whose graded Betti numbers are an upper bound for those of the monomial ideal $I$.  We prove that strongly stable ideals are fixed under these operations.  Section 3 contains the proof of our main theorem using a result of Caviglia and Kummini \cite{CK} to reduce the problem to the characteristic zero case.

\section{Shifting and Compression}

Throughout this paper, $A=K[x_1,...,x_n]$ is a polynomial ring over a field $K$, where char($K$) is arbitrary, and $P=(x_1^{e_1},...,x_r^{e_r})$, for some $r\leq n$ and $2\leq e_1\leq e_2\leq \cdots \leq e_r$ .  Furthermore, throughout this section, assume $I$ is a monomial ideal containing $P+J$ where $J$ is a strongly stable ideal.  Recall that $J$ is strongly stable if it satisfies the combinatorial property that whenever $x_i m\in J$, then $x_j m\in J$ for all monomials $m$ and for all $j<i$ \cite{G}.

In the proof of their main theorem, Mermin and Murai show that there exists a strongly-stable-plus-$P$ ideal $B$ with the same Hilbert function as the ideal $I$ such that $b_{ij}(B)\geq b_{ij}(I)$ for all $i,j$,  \cite[Proposition 8.7]{MM}.  In this section, we recall the operations that Mermin and Murai use to construct the ideal $B$ from $I$, and show that in addition to the properties above, we also have $J\subset B$.  We will conclude this section with the proof of the following proposition:

\begin{Proposition}\label{MainProp}
If $I$ is a monomial ideal containing $P+J$, then there exists a strongly-stable-plus-$P$ ideal $B$ with the same Hilbert function as $I$ such that $b_{ij}(B)\geq b_{ij}(I)$ for all $i,j$ and $J\subset B$.
\end{Proposition}

For pairs of variables $a>_{\text{lex}} b$, the ideal $B$ is constructed in \cite{MM} in finitely many steps by replacing $I$ with any of the following ideals:
\begin{enumerate}
\item Shift$_{a,b}(I)$,
\item Shift$_{a,b,t}(I)+P$,
\item $T=T'+P$ as in Proposition \ref{Compression}.
\end{enumerate}

We introduce the definitions of the basic operations used above, and prove that strongly stable ideals do not move after replacing $I$ by any of these ideals.

\begin{Definition}
Let $I$ be a monomial ideal, and fix variables $a>_{lex} b$ and $t\in \mathbb{Z}_{\geq 0}$. The {\bf $(a,b,t)$-shift} of $I$, denoted Shift$_{a,b,t}(I)$, is the $K$-vector space generated by monomials of the form:
\begin{equation*}
\begin{Bmatrix*}
\begin{split}
fa^sb^r &| fa^sb^r\in I, r<t\\
fa^sb^{s+t} &| fa^sb^{s+t}\in I\\
fa^lb^{s+t} &| fa^lb^{s+t} \in I \text{ or } fa^sb^{l+t}\in I\\
fa^sb^{l+t} &| fa^lb^{s+t} \in I \text{ and } fa^sb^{l+t}\in I
\end{split}
\end{Bmatrix*}
\end{equation*}
where the set is taken over all monomials $f$ such that $a\nmid f$ and $b\nmid f$, and over all integers $0\leq s<l$.
\end{Definition}

\begin{Remark}
Notice that when $fa^lb^{s+t} \in I$ and $fa^sb^{l+t}\in I$, then both monomials $fa^lb^{s+t}$ and $fa^sb^{l+t}$ will be generators of Shift$_{a,b,t}(I)$.
\end{Remark}

\begin{Definition}
The $(a,b)$-shift of $I$ is the $(a,b,0)$-shift of $I$ as defined above.
\end{Definition}

\begin{Remark}  For $t\neq 0$, Shift$_{a,b,t}(I)$ does not necessarily fix ideals generated by powers of variables.  Thus, in order to preserve the ideal $P$ when applying the shifting operation for $t\neq 0$, Mermin and Murai use the operation Shift$_{a,b,t}(I)+P$.
\end{Remark}

\begin{Proposition}\label{Shifting}
Let $I$ be a monomial ideal containing $P+J$.  Fix variables $a >_{\text{lex}} b$ and $t>0$.  Then, $J\subset \text{Shift}_{a,b,t}(I)$.
\end{Proposition}

\textit{Proof.} Write $m=m'a^{\alpha}b^{\beta}\in J$, where $a\nmid m'$ and $b\nmid m'$.  If $\beta\leq \alpha+t$, then it's clear that $m\in \text{Shift}_{a,b,t}(I)$.  The only case where we need to use the assumption that $J$ is strongly stable is when $\beta >\alpha+t$.  Here, we need to show that $m'a^{\beta-t}b^{\alpha+t}\in I$.  Let $N=\beta-(\alpha+t)$.  Since $J$ is strongly stable and $N>0$, then $m\cdot \frac{a^N}{b^N}\in J\subset I$.  We see that $m\cdot\frac{a^N}{b^N}=m'a^{\beta-t}b^{\alpha+t}$.  Since both $m=m'a^{\alpha}b^{(\beta-t)+t}\in I$ and $m'a^{\beta-t}b^{\alpha+t}\in I$, it follows that $m\in\text{Shift}_{a,b,t}(I)$. \qed

The final operation used to transform the ideal $I$ in the proof of Mermin and Murai is a compression.  The following definition is described by Mermin in \cite{Me}:

\begin{Definition}
Let $I$ be a monomial ideal, and fix variables $a>_{lex} b$.  Write $I$ as a direct sum of the form $I=\bigoplus\limits_f fV_f$, where the sum is taken over all monomials $f$ in $K[\{x_1,...,x_n\}\setminus \{a,b\}]$ and $V_f$ are $K[a,b]$-ideals.  The {\bf $\{a,b\}$-compression} of $I$ is the ideal $\bigoplus\limits_f fN_f$, where $N_f\subset K[a,b]$ are the lex ideals with the same Hilbert function as $V_f$.
\end{Definition}

\begin{Proposition}\label{Compression}
Let $I$ be a monomial ideal containing $P+J$.  Fix variables $a >_{\text{lex}} b$.  Let $I'$ be the ideal of $A$ generated by all the minimal generators of $I$ except for $b^{e_b}$, let $T'$ be the $\{a,b\}$-compression of $I'$, and let $T=T'+P$.  Then, $J\subset T$.
\end{Proposition}

\textit{Proof.}  As in the definition of $\{a,b\}$-compression, write $I'=\bigoplus\limits_f fV_f$ with $f\in \text{Mon}(K[\{x_1,...,x_n\}\setminus \{a,b\}])$ and $V_f\subset K[a,b]$.  Let $T'=\bigoplus\limits_f fN_f$ be the $\{a,b\}$-compression of $I'$.  First, suppose $b^{e_b}$ is not a minimal generator of $I$.  In this case, $I'=I$, and therefore, $T'$ is the $\{a,b\}$-compression of $I$.  Since strongly stable ideals are $\{a,b\}$-compressed, as stated in Proposition 3.8 of \cite{Me}, then $J\subset T'$.

If instead $b^{e_b}$ is a minimal generator of $I$, let $m=m'a^{\alpha}b^{\beta}$ be a monomial in $J$ with $a\nmid m'$, $b\nmid m'$.  Clearly, if $\beta \geq e_b$, then $m\in P\subset T$.  So we may assume $\beta <e_b$.  Since $J$ is strongly stable, then we have:
\begin{equation*}
m=m'a^{\alpha}b^{\beta}<_{lex} m'a^{\alpha +1}b^{\beta -1}<_{lex} \cdots <_{lex} m'a^{\alpha +\beta}\in J.
\end{equation*}
Furthermore, all of these monomials are in $I'$.  Thus,
\begin{equation*}
a^{\alpha}b^{\beta}<_{lex}a^{\alpha +1}b^{\beta -1}<_{lex} \cdots <_{lex}a^{\alpha +\beta}\in V_{m'}.
\end{equation*}
These are the first monomials of degree $\alpha +\beta$ in $K[a,b]$, hence they are also elements of the lex ideal $N_{m'}$.  In particular, this implies that $m\in T'$. \qed\\

We conclude this section with the proof of the main propostion:\\

\textit{Proof of Proposition \ref{MainProp}.} By Proposition 8.7 of \cite{MM}, there exists a strongly-stable-plus-$P$ ideal $B$ with the same Hilbert function as $I$ and $b_{ij}(B)\geq b_{ij}(I)$ for all $i,j$. Furthermore, by Propositions \ref{Shifting} and \ref{Compression}, strongly stable ideals do not move under the operations used to construct the ideal $B$.  Hence, when $J\subset I$, we also have $J\subset B$.\qed

\section{Main Result}

In the previous section, we showed that strongly stable ideals do not move under any of the three above operations.  Now, we apply this to the situation in which $J=\widetilde{L}$ is a piecewise lex ideal.  We begin by reminding the reader of the definition introduced by Shakin \cite{S}.

\begin{Definition}
For each $1\leq i\leq n$, let $A_{(i)}$ be the polynomial ring over $K$ in the first $i$ variables.  An ideal $\widetilde{L}\subset A$ is called a {\bf piecewise lex ideal} if it can be written as a sum:
\begin{equation*}
\widetilde{L}=L_{(1)}A+L_{(2)}A+...+L_{(n)}A
\end{equation*}
where $L_{(i)}$ is a lex ideal in the ring $A_{(i)}$ for each $i$.
\end{Definition}

Since piecewise lex ideals are strongly stable, then we have shown that when $I$ contains $P+\widetilde{L}$, the ideal $B$ does as well.

\begin{Theorem}
Let $I \subset A$ be a homogeneous ideal with $P+\widetilde{L}\subset I$.  There exists a lex ideal $L$ such that
\begin{enumerate}[(i)]
\item $P+\widetilde{L}+L$ has the same Hilbert function as $I$.
\item $b_{ij}(P+\widetilde{L}+L)\geq b_{ij}(I)$ for all $i, j$
\end{enumerate}
\end{Theorem}

\textit{Proof.} Without loss of generality, using a standard upper-semicontinuity argument, we may assume $I$ is a monomial ideal containing $P+\widetilde{L}$ by replacing $I$ with in$(I)$.  By Proposition \ref{MainProp}, there is a strongly-stable-plus-$P$ ideal $B$ with the same Hilbert function as $I$ and such that $b_{ij}(B)\geq b_{ij}(I)$.  Furthermore, we have that $P+\tilde{L}\subset B$. Since $B$ is a strongly-stable-plus-$P$ ideal, the graded Betti numbers, $b_{ij}(B)$, do not depend on char($K$) by \cite[Corollary 3.7]{CK}.  Hence, we can assume char($K$)=0.  The characteristic zero result of Caviglia and Sbarra, \cite[Theorem 3.4]{CS}, gives a lex ideal $L$ such that $P+\widetilde{L}+L$ has the same Hilbert function as $B$ and $b_{ij}(P+\widetilde{L}+L)\geq b_{ij}(B)$ for all $i, j$.  Again, since $P+\widetilde{L}+L$ is strongly-stable-plus-$P$, then the Betti numbers do not depend on the characteristic, so the inequality also holds for char($K$) arbitrary. \qed


\begin{bibdiv}
\begin{biblist}

\bib{A}{article}{
	title={On the Eisenbud-Green-Harris conjecture}
	author={Abedelfatah, A.}
	date={2015}
	journal={Proc. Amer. Math. Soc.}
	volume={143}
	pages={105--115}
}

\bib{B}{article}{
	title={Upper bounds for the Betti numbers of a given Hilbert function}
	author={Bigatti, A.}
	date={1993}
	journal={Comm. Algebra}
	volume={21}
	pages={2317--2334}
}

\bib{CK}{article}{
     title={Poset embeddings of Hilbert functions and Betti numbers}
     author={Caviglia, G.}
     author={Kummini, M.}
     date={2014}
     journal={J. Algebra}
     volume={410}
     pages={244--257}
}

\bib{CM}{article}{
	title={Some cases of the Eisenbud-Green-Harris conjecture}
	author={Caviglia, G.}
	author={Maclagan, D.}
	date={2008}
	journal={Math. Res. Lett.}
	volume={15}
	pages={427--433}
}

\bib{CS}{article}{
     title={Distractions of Shakin rings}
     author={Caviglia, G.}
     author={Sbarra, E.}
     date={2014}
     journal={J. Algebra}
     volume={419}
     pages={318--331}
}

\bib{C}{article}{
     title={An application of liaison theory to the Eisenbud-Green-Harris conjecture}
     author={Chong, K. F. E.}
     note={arXiv:1311.0939v1 [math.AC]}
}

\bib{CL}{article}{
     title={A generalization of a combinatorial theorem of Macaulay}
     author={Clements, G.F.}
     author={Lindstr{\"o}m, B.}
     date={1969}
     journal={J. Combinatorial Theory}
     volume={7}
     pages={230--238}
}

\bib{EGH1}{article}{
	title={Higher Castelnuovo theory}
	author={Eisenbud, D.}
	author={Green, M.}
	author={Harris, J.}
	date={1993}
	journal={Ast\`e risque}
	volume={218}
	pages={187--202}
}

\bib{EGH2}{article}{
	title={Cayley-Bacharach theorems and conjectures}
	author={Eisenbud, D.}
	author={Green, M.}
	author={Harris, J.}
	date={1996}
	journal={Bull. Amer. Math. Soc}
	volume={33}
	pages={295--324}
}

\bib{F}{article}{
	title={Almost complete intersections and the lex-plus-powers conjecture}
	author={Francisco, C.}
	date={2004}
	journal={J. Algebra}
	volume={276}
	pages={737--760}
}

\bib{FR}{article}{
	title={Lex-plus-powers ideals}
	author={Francisco, C.}
	author={Richert, B.}
	date={2007}
	journal={Lect. Notes Pure Appl. Math.}
	volume={254}
	pages={113--144}
}

\bib{G}{article}{
	title={Generic initial ideals}
	author={Green, M.}
	date={2010}
	booktitle={Six lectures on commutative algebra}
	series={Mod. Birkh{\"a}user Class.}
	publisher={Birkh{\"a}user Boston}
	address={Boston}
	pages={119--186}
}

\bib{H}{article}{
	title={Maximum Betti numbers of homogeneous ideals with a given Hilbert function}
	author={Hulett, H.}
	date={1993}
	journal={Comm. Algebra}
	volume={21}
	pages={2335--2350}
}

\bib{K}{article}{
	title={Algebraic shifting}
	author={Kalai, G.}
	date={2002}
	journal={Adv. Stud. Pure Math.}
	volume={33}
	pages={121--163}
}

\bib{Ma}{article}{
    label={Ma}
	title={Some properties of enumeration in the theory of modular systems}
	author={Macaulay, F.}
	date={1927}
	journal={Proc. London Math. Soc.}
	volume={26}
	pages={531--555}
}

\bib{Me}{article}{
    label={Me}
    title={Compressed ideals}
    author={Mermin, J.}
    date={2008}
    journal={Bull. London Math. Soc.}
    volume={40}
    pages={77--87}
}

\bib{MM}{article}{
     title={The Lex-Plus-Powers Conjecture holds for pure powers}
     author={Mermin, J.}
     author={Murai, S.}
     date={2011}
     journal={Adv. Math.}
     volume={226}
     number={3}
     pages={3511--3539}
}

\bib{MP}{article}{
	title={Lexifying ideals}
	author={Mermin, J.}
	author={Peeva, I.}
	date={2006}
	journal={Math. Res. Lett.}
	volume={13}
	pages={409--422}
}

\bib{MH}{article}{
	title={Algebraic Shifting and graded Betti numbers}
	author={Murai, S.}
	author={Hibi, T.}
	date={2009}
	journal={Trans. Amer. Math. Soc.}
	volume={361}
	pages={1853--1865}
}

\bib{P}{article}{
	title={Deformation classes of graded modules and maximal Betti numbers}
	author={Pardue, K.}
	date={1996}
	journal={Illinois J. Math.}
	volume={40}
	pages={564--585}
}

\bib{R}{article}{
	title={A study of the lex plus powers conjecture}
	author={Richert, B.}
	date={2004}
	journal={J. Pure Appl. Algebra}
	volume={186}
	pages={169--183}
}

\bib{S}{article}{
	title={Piecewise lexsegment ideals}
	author={Shakin, D.}
	date={2003}
	journal={Sb. Math.}
	volume={194}
	pages={1701--1724}
}

\end{biblist}
\end{bibdiv}
\end{document}